\documentclass[12pt]{article}

\usepackage{amsfonts}

\newtheorem{thm}{Theorem}[section]
\newtheorem{lemma}[thm]{Lemma}
\newtheorem{cor}[thm]{Corollary}

\newtheorem{prop}[thm]{Proposition}

\newcommand{\proof
}{\par\medskip\noindent {\bf Proof.\ \ }}

\newcommand{\be}{\begin{equation}}
\newcommand{\ee}{\end{equation}}
\newcommand{\openbox}{\leavevmode
  \hbox to8pt{\hfil\vrule\vbox to6pt{\hrule width6pt\vfil\hrule}\vrule}}

\newcommand{\qed}{\hbox to5pt{ } \hfill \openbox\bigskip\medskip}

\newcommand{\cX}{\mbox{$\cal X$}}

\newcommand{\cF}{\mbox{$\cal F$}}
\newcommand{\cG}{\mbox{$\cal G$}}
\newcommand{\cH}{\mbox{$\cal H$}}

\newcommand{\cM}{\mbox{$\cal M$}}
\newcommand{\cN}{\mbox{$\cal N$}}
\newcommand{\cV}{\mbox{$\cal V$}}
\newcommand{\cB}{\mbox{$\cal B$}}

\newcommand{\cR}{\mbox{$\cal R$}}
\newcommand{\cU}{\mbox{$\cal U$}}
\newcommand{\cW}{\mbox{$\cal W$}}

\newcommand{\Q}{\mathbb Q}

\newcommand{\F}{\mathbb F}

\newcommand{\ve}[1]{\mathbf{#1}}
\newcommand{\monom}[2]{\ve{#1}^{\ve{#2}}}

\title{Multivalued generalizations 
of the Frankl--Pach Theorem}
\author{G\'abor Heged\H{u}s
\\ {\small 
Johann Radon Institute for Computational and Applied Mathematics, Linz}
\medskip
\\ Lajos R\'onyai
\\ {\small Computer and Automation Institute, Hungarian Academy of Sciences}
\\ {\small and }
\\ {\small Institute of Math., Budapest University of Technology and
Economics}
}

\begin{document}

\footnotetext{Research supported in part by OTKA grants NK 63066, NK
72845, K77476, and K77778.}
\maketitle

\begin{abstract}
In \cite{FP} P. Frankl and J. Pach proved the following uniform version of 
Sauer's Lemma.

\medskip
\noindent
{\em Let $n,d,s$ be natural numbers such that $d\leq n$, $s+1\leq n/2$. 
Let $\cF
\subseteq {[n] \choose d}$ be an arbitrary $d$-uniform set 
system such that
$\cF$ does not shatter an $s+1$-element set, then
$$
|\cF|\leq {n \choose s}.
$$  }
\medskip

We prove here two generalizations of the above theorem to 
$n$-tuple systems. To obtain these results, we use Gr\"obner basis methods, 
and describe the standard monomials of Hamming spheres.

\end{abstract}


\medskip

\section{Introduction}

Let  $[n]$ stand for the set $\{1,2,
\ldots, n\}$. 
The family of all subsets of $[n]$ is denoted by $2^{[n]}$. 
For an integer $0\leq d\leq n$ we denote by
${[n] \choose d}$ the family of all  $d$ element subsets of $[n]$, and
by ${[n] \choose \leq d}={[n] \choose 0}\cup\cdots\cup{[n] \choose d}$
the family of subsets of size at most $d$.

Let $n>0$, $\cF\subseteq 2^{[n]}$ be a family of subsets of $[n]$, 
and $S$ be a subset of $[n]$. We say that $\cF$ {\it shatters} $S$ if
\begin{equation} \label{sha1}
\{F\cap S :~F\in \cF\}=2^{S}.
\end{equation}
Define
\begin{equation}
\mbox{sh} (\cF)=\{S\subseteq [n] :~\cF\mbox{\ shatters\ }S\}.
\end{equation}

The following result was proved by 
Sauer, \cite{S}, and independently by  Vapnik and Chervonenkis \cite{VC},
and Perles and Shelah \cite{Sh}:

\begin{thm} \label{VC}
Suppose that $0\leq s \leq n-1$ and let $\cF\subseteq 2^{[n]}$ be 
an arbitrary set family with no shattered set of size $s+1$. Then
$$
|\cF|\leq \sum_{i=0}^{s} {n \choose i}.
$$
\end{thm}

\medskip

Karpovsky and Milman in \cite{KM} gave a generalization of 
Sauer's result for tuple systems. Next we explain this multivalued
generalization.  Throughout the paper $q\geq 2$ is an integer.
Let  $(q)$ stand for the set $\{0,1,
\ldots, q-1\}$ and denote by $\ve v_F$ the characteristic vector of the set 
$F\subseteq [n]$. Clearly we have $\ve v_F\in (2)^n$. 

Subsets  $\cV\subseteq (q)^n$ will be called {\em tuple
systems}\footnote{They are also called {\em sets of vectors} in the
literature.}. Note 
that an element $\ve v$ of a tuple system $\cV$ can also be viewed as a 
function
from $[n]$ to $(q)$. With this in mind,  we say  
that the tuple  system $\cV$ shatters a set $S\subseteq [n]$, if 
$$
\{ \ve v\mid_{S}:~ \ve v\in \cV\}
$$
is the set of all functions from $S$ to $(q)$; here $\ve v\mid_{S}$ denotes the 
restriction of the function  $\ve v$ to the set $S$. This extends the
binary notion of shattering introduced in (\ref{sha1}). 
In fact, consider
\begin{equation}
\mbox{Sh} (\cV):=\{S\subseteq [n]:~ \cV\mbox{\ shatters\ }S\},
\end{equation}
the set of the shattered sets of the tuple system $\cV$. 

Clearly $\mbox{Sh}(\cV)\subseteq 2^{[n]}$.
Moreover, if $\cF\subseteq 2^{[n]}$ is a set system, then 
$$
\mbox{sh}(\cF)= \mbox{Sh}(\{\ve v_F\in 2^{(n)}:~  F\in \cF\}).
$$

The following result was proved by Karpovsky and Milman in 
\cite[Theorem 2]{KM} (see also Alon
\cite[Corollary 1]{A}, Steel \cite[Theorem 2.1]{S} and 
Anstee \cite[Theorem 1.3]{An}).

\begin{thm} \label{KaMi}
Let $0\leq s\leq n-1$ be an integer and let $\cV\subseteq (q)^n$ be a
tuple system with no shattered set of size $s+1$. 
Then
$$
|\cV|\leq \sum_{i=0}^s (q-1)^{n-i} {n \choose i}.
$$ 
\end{thm} \qed

The above theorem can be viewed as a natural multivalued generalization of 
Theorem \ref{VC}. 

A set family ${\cal F}\subseteq 2^{[n]}$ is called $d$-uniform, iff $|F|=d$
holds, whenever $F\in {\cal F}$. Uniformity can be generalized to tuple
systems in two simple ways.  First let $0\leq d\leq (q-1)n$. A tuple system
$\cV\subseteq (q)^n$ is {\em $d$-uniform}  iff $\sum_{i=1}^n v_i=d$ holds for 
every 
$(v_1,\ldots,v_n)\in \cV$. 

Alternatively, let $0\leq d\leq n$. A tuple system $\cV \subseteq (q)^n$ 
is {\em $d$-Hamming}, if $ |\{i\in[n]:~ v_i\neq 0\}|=d$ for every 
$(v_1,\ldots,v_n)\in \cV$. 

In \cite{FP} P. Frankl and J. Pach proved the following uniform version of
Theorem \ref{VC}.

\begin{thm}
Let $n,d,s$ be natural numbers such that $d\leq n$, $s+1\leq n/2$. Let $\cF
\subseteq {[n] \choose d}$ be an arbitrary $d$-uniform set system
such that
$\cF$ does not shatter an $s+1$-element set, then
$$
|\cF|\leq {n \choose s}.
$$
\end{thm}

We would like to extend this result to tuple systems and hence obtain 
uniform variants of the Karpovsky--Milman theorem. We prove the following 
two theorems, which specialize to the Frankl--Pach bound in the case 
$q=2$. 

\begin{thm} \label{main1}
Suppose that $0\leq d\leq n(q-1)$ and $s\leq \frac{n}{2}$. Let 
$\cV$ be an arbitrary $d$-uniform tuple system with no shattered set of 
size $s+1$.
Then
$$
|\cV|\leq \sum_{i=0}^s (q-1)^{n-i}\left( {n \choose i}-{n \choose i-1}\right) .
$$
\end{thm}

\medskip

\begin{thm} \label{main2}
Suppose that $0\leq d\leq n$ and $0\leq d+s\leq n$. Let 
$\cV$ be an arbitrary $d$--Hamming tuple system 
with no shattered set of size $s+1$.
Then
$$
|\cV|\leq {n\choose s} \sum_{i=0}^d {n-s \choose i}(q-2)^i.
$$
\end{thm}

The rest of the paper is organized as follows: Section 2 contains our basic
results involving Gr\"obner bases and normal sets. Sections 3 and 4 contain
the proofs of Theorems \ref{main1} and \ref{main2}. The paper ends with some
concluding remarks.

\medskip

{\bf Acknowledgements.} We thank the referees for their helpful comments.

\section{Gr\"obner bases, standard monomials and shattering }

Next we fix some notation related to  Gr\"obner 
bases in polynomial rings, we need later on. The interested reader 
can find a detailed introduction to this topic in the
classic papers by Buchberger \cite{B1}, \cite{B2}, \cite{B3},
and in the textbooks \cite{AL}, \cite{CCS}, \cite{CLO}.

We shall work over the field of rational numbers $\Q$ and 
we denote by  $R:=\Q[x_1,\ldots ,x_n]$ the polynomial ring in $n$ variables over 
$\Q$. We fix a monomial order $\prec$ on $R$ such that $x_n\prec
x_{n-1}\prec \cdots \prec x_1$ holds. For a nonzero polynomial $f\in R$ we
denote by  ${\rm lm}(f)$ the largest monomial of $f$ with respect to
$\prec$. 

Let $I$ be a nonzero ideal of $R$. Recall, that a finite 
subset $\cG\subseteq I$ is a {\it
Gr\"obner basis} of $I$ (with respect to $\prec$) if for 
every $f\in I$ there exists a $g\in \cG$ such
that ${\rm lm}(g)$ divides ${\rm lm}(f)$. 

We shall denote by  ${\rm SM}(I)$ the set of all standard monomials of $I$ with 
respect to the term-order $\prec$ over $\Q$. ${\rm SM}(I)$ is 
often called as a {\em normal set} of $I$. ${\rm SM}(I)$ is the complement
of ${\rm LM}(I)$, the set of all leading monomials for $I$ within the set of
all monomials of $R$.  
It is known  that for a nonzero ideal $I$  (the image of) 
${\rm SM}(I)$ is a basis of the $\Q$-vector-space $R/I$. 

We denote by $NF(f,\cG)$ the (unique) normal form of a polynomial $f\in R$ 
with respect to a Gr\"obner basis $\cG$.


To study the polynomial functions on a (finite) set of vectors 
$\cV\subseteq \Q^n$,
it is convenient to work with the ideal $I(\cV)$:
$$ I(\cV):=\{f\in R:~f(\ve v)=0 \mbox{ whenever } \ve v \in \cV\}. $$

It is immediate that $\mbox{SM}(I(\cV))$ is downward closed: if 
$y\in {\rm
SM}(I(\cV))$, $y_1,y_2$ are
monomials from $R$ such that $y=y_1y_2$ then 
$y_1\in {\rm SM}(I(\cV))$.

An easy interpolation argument shows that any function from $\cV$ to $\Q$ 
is a polynomial. This gives a bijection from $\cV $ to  ${\rm SM}(I(\cV))$. 
We obtain in particular, that 
\begin{equation} \label{egyenlo}
|\mbox{SM}(I(\cV))|=|\cV|.
\end{equation}

\subsection{Standard monomials and shattering}

The following example illustrates some of the notions we have mentioned so
far. Also, it will be useful later in the paper. 

\medskip

\noindent
{\bf Example.} We describe a Gr\"obner
basis and the standard monomials of the set $(q)^n\subseteq \Q^n$. 
We introduce the polynomials
\begin{equation} \label{q-poli}
f_i(x_i):= \prod_{j=0}^{q-1} (x_i-j)\in \Q[x_1,\ldots ,x_n]
\end{equation}
for $1\leq i\leq n$. These polynomials vanish on $(q)^n$, and the leading 
monomial of  $f_i(x_i)$ is $x_i^q$. These imply that $\mbox{SM}(I ((q)^n))$ is
a subset of $\{\monom xv:~ \ve v\in (q)^n\}$. But this latter set 
has $q^n$ elements,
hence by (\ref{egyenlo}) we have 
\begin{equation} \label{q-standard}
\mbox{SM}(I( (q)^n))=\{\monom xv:~\ve v\in (q)^n\}. 
\end{equation}
This in turn implies that $\cG= \{f_1(x_1), \ldots ,f_n(x_n)\}$ is a 
Gr\"obner basis  for $I((q)^n)$. 

\medskip

Next we prove a statement, which connects the notion of shattering 
to the theory of Gr\"obner bases.

\begin{prop} \label{alap}
Let $\cV\subseteq (q)^n$ be a set of tuples. If 
$S=\{i_1,\ldots, i_k\}\subseteq [n]$ is a set for which 
$x_{i_1}^{q-1}\cdots x_{i_k}^{q-1} \in 
\mbox{\rm SM}(I(\cV))$, 
then $S\in \mbox{\rm Sh}(\cV)$. 
\end{prop}

\proof  Suppose that $S\notin \mbox{Sh}(\cV)$. 
We  show that $x_{i_1}^{q-1}\cdots x_{i_k}^{q-1} \notin 
\mbox{SM}(I(\cV))$. 
As $S\notin \mbox{Sh}(\cV)$, there exists a tuple 
$\ve w=(w_1,\ldots,w_n)\in (q)^n$
such that  $\ve w\mid_{S}\neq \ve v\mid_{S}$ holds for every $\ve v\in \cV$. 


Consider now the polynomial
$$
g(x_1,\ldots,x_n):=\prod_{j\in S} h_j(x_j)\in \Q[x_1,\ldots,x_n],
$$
where
$$
h_j(x_j):=\prod_{i=0,i\neq w_j}^{q-1} (x_j-i)\in \Q[x_j].
$$
Then we immediately see that  
\begin{equation} \label{fotag}
\mbox{lm}(g)=x_{i_1}^{q-1}\cdots x_{i_k}^{q-1}.
\end{equation}
We claim that $g(\ve v)=0$ holds for every  $\ve v\in \cV$.
Indeed, let $\ve v=(v_1,\ldots,v_n)\in\cV$ be an arbitrary tuple. Since
$\ve w\mid_{S}\neq \ve v\mid_{S}$, there must exist an index $j\in S$ such 
that $w_j\neq v_j$.
Then 
\begin{equation}
h_j(v_j)=\prod_{i=0,i\neq w_j}^{q-1} (v_j-i)=0,
\end{equation}
implying that $g(\ve v)=0$. We obtained that
$g\in I(\cV)$. This, together with 
(\ref{fotag}) implies that 
$$
x_{i_1}^{q-1}\cdots x_{i_k}^{q-1}=
\mbox{lm}(g)
\notin 
\mbox{Sm}(I(\cV)).
$$
\qed

\subsection{The blow-up of a set family}

Let $\ve v\in \Q^n$ be an $n$-tuple, and put 
$$
\mbox{supp}(\ve v):=\{i\in [n]:~ v_i\neq 0\}.
$$

Let $\cF\subseteq 2^{[n]}$ be a set system. We define 
the {\em blow--up} $\cF^q\subseteq (q)^n$ 
of $\cF$ as
$$
\cF^q:= \{\ve v\in (q)^n:~ \mbox{supp}(\ve v)\in \cF\}.
$$

Clearly 
$$
|\cF^q|=\sum_{F\in \cF} (q-1)^{|F|}.
$$

For a subset $J\subseteq [n]$, we consider
$$
\cF_J:= \{F\in \cF:~ J\subseteq F \}\subseteq  2^{[n]}.
$$
Let $g(x_1,\ldots ,x_n)\in \Q[x_1,\ldots ,x_n]$ be a polynomial.
We define $\overline{g}(x_1,\ldots ,x_n):=g(p(x_1),\ldots ,p(x_n))$, 
where $p(x)\in
\Q[x]$ is the unique polynomial  for which  $\mbox{deg}(p)=q-1$, 
$p(0)=0$ and $p(i)=1$ for each $1\leq i\leq q-1$. 
Clearly we have $\mbox{lm}(\overline{g})=\mbox{lm}(g)^{q-1}$.

For a tuple $\ve v=(v_1,\ldots ,v_n) \in (q)^n$ we define three 
subsets $J(\ve v),
Q(\ve v)$ and $Z(\ve v)$ of $[n]$ as follows:
$$ J(\ve v)=\{ i\in [n]:~0<v_i<q-1 \},~~~
Q(\ve v)=\{ i\in [n]:~v_i=q-1 \}, \mbox{ and } $$
$Z(\ve v)=\{i\in [n]:~v_i=0 \}$. The sets $J(\ve v),Q(\ve v),Z(\ve v)$ 
partition $[n]$.
We note also that a set family $\cF\subseteq 2^{[n]}$ can be 
identified with the tuple system 
$$ \{\ve v_F:~F \in \cF \} \subseteq (2)^n. $$ 
Here $\ve v_F$ denotes the characteristic vector of a set 
$F\subseteq [n]$. This way we can speak
of Gr\"obner bases and standard monomials for a set family $\cF$. 
 
The next result, which may be of independent interest, relates the Gr\"obner
bases and normal sets of $\cF^q$ to those of the set systems $\cF_J$,
$J\subseteq [n]$.  It establishes a useful connection of the multivalued 
case to the sometimes simpler binary case. We recall first that the
polynomials $f_1, \ldots ,f_n$ from (\ref{q-poli}) form a Gr\"obner basis
of the ideal of $(q)^n$. 

For a subset $J\subseteq [n]$ $x_J$ denotes the monomial $x_J:=\prod _{j\in
J}x_j$. In particular, $x_{\emptyset}=1$.

\begin{thm} \label{monom}
Let $\cF\subseteq [n]$ be a nonempty set family. 
For $J\subseteq [n]$, let $\cG(\cF_J)$ denote a fixed Gr\"obner basis of the 
ideal $I(V(\cF_J))$. Then
\begin{equation} \label{unio2}
\{f_1,\ldots ,f_n\}\cup (\cup_{J\subseteq [n]} 
\{x_J\cdot \overline{g}:~ g\in \cG(\cF_J)\})\cup \{x_J:~ J\subseteq[n],~ 
\cF_J=\emptyset \}
\end{equation}
is a Gr\"obner basis of the ideal $I(\cF^q)$.
Moreover,
\begin{equation} \label{stan} 
{\rm SM}(I(\cF^q))=\{\monom xv:~\ve v\in (q)^n,~~ 
\cF_{J(\ve v)} \neq \emptyset,\ ~{\rm and }~x_{Q(\ve v)} \in {\rm
SM}(I(\cF_{J(\ve v)}))
\}. 
\end{equation}
\end{thm}

\proof
We note first that the polynomials from (\ref{unio2}) clearly 
vanish on $\cF^q$.
Let $\cR$ denote the right hand side of (\ref{stan}). To establish the
Theorem, it suffices to prove
that $|\cR|=|\cF^q|$, and for each $y=\monom xv \notin \cR$ 
there exists a polynomial $h$ from the set (\ref{unio2})
such that the leading monomial of $h$ divides $y$.

Indeed, then $y\in \mbox{LM}(I(\cF^q))$. Using also
(\ref{q-standard}) we
obtain that $\mbox{SM}(I(\cF^q))\subseteq \cR$. 
But then  $|\cR|=|\cF^q|=|\mbox{SM}(I(\cF^q))|$ implies 
that $\mbox{SM}(I(\cF^q))=\cR$ and in turn gives that the 
union (\ref{unio2}) constitutes a Gr\"obner basis of the ideal $I(\cF^q)$.

\medskip

First we prove that $|\cR|=|\cF^q|$.
For each $J\subseteq [n]$ such that $\cF_J\not=\emptyset $
we fix a bijection 
$$
\phi_J:\{F\in \cF:~ J\subseteq F\}\to \mbox{SM}(I( \cF_J)).
$$ 
From (\ref{egyenlo}) we see that  
$$|\mbox{SM}(I(\cF_J))|=|\cF_J|=|\{F\in
\cF:~ J\subseteq F\}|,$$
hence such maps exist. Next we show that  
the following is a disjoint union decomposition of $\cR$:
\begin{equation} \label{unio3}
\cR=\cup_{F\in \cF} \cup_{J\subseteq F}\ \{\monom xv :~ \ve v\in
(q)^n,~~J(\ve v)=J
\mbox{ and }x_{Q(\ve v)}=\phi_J(F)\}.
\end{equation}
Indeed, a monomial $\monom xv$ from the right side belongs to $\cR$, 
because $\phi _J(F)$ is in   
${\rm SM}(I(\cF_{J(\ve v)}))$.  Conversely, if 
$\monom xv \in \cR$, then 
$x_{Q(\ve v)}=\phi _J(F)$ for some $F\in \cF $ with $J\subseteq F$, because 
$\phi _J$ is surjective. 

Let $J\subseteq F \subseteq [n]$ be fixed subsets, with $F\in \cF$. 
Then $\cF_J\neq \emptyset$ and we have 
\begin{equation} \label{runio}
|\{\monom xv \in \cR:~  J(\ve v)=
J\mbox{ and }x_{Q(\ve v)}=\phi_J(F)\}|=(q-2)^{|J|}.
\end{equation}
Keeping this in mind, for a fixed $F\in \cF$ we have
$$
|\cup_{J\subseteq F} \{\monom xv :~~\ve v\in (q)^n, ~~J(\ve v)=J
\mbox{ and } x_{Q(\ve v)}=\phi_J(F) \}|=
$$
$$
\sum_{J\subseteq F} |\{\monom xv :~~\ve v\in (q)^n, ~~J(\ve v)=J
\mbox{ and } x_{Q(\ve v)}=\phi_J(F) \}|=
$$
$$
\sum_{J\subseteq F} (q-2)^{|J|}=\sum_{i=0}^{|F|} {|F| 
\choose i}(q-2)^i =(q-1)^{|F|}.
$$
Using again that (\ref{unio3}) is a disjoint decomposition, we infer that 
$$
|\cR|=\sum_{F\in \cF} (q-1)^{|F|}=|\cF^q|.
$$

Finally, we prove that if $y=\monom xv \notin \cR$, then 
$y\in \mbox{LM}(I(\cF^q))$, more precisely, 
$y$ is divided by the leading monomial of some polynomial $h$ from 
(\ref{unio2}).

If $v_i>q-1$, then $h=f_i(x_i)$ will do. We can therefore assume, 
that $\ve v\in
(q)^n$. Now if $\cF_{J(\ve v)}=\emptyset$, then $h=x_{J(\ve v)}$ is 
a good choice.
We are left with the case  $\cF_{J(\ve v)} \not= \emptyset$. Then
$\monom xv \not\in\cR$ is possible only if $x_{Q(\ve v)}$ is a 
leading monomial for
the ideal $I(\cF_{J(\ve v)})$, hence there exists a $g\in \cG (\cF_{J(\ve v)})$ 
whose leading term divides  $x_{Q(\ve v)}$. Taking also into consideration 
that $x_{J(\ve v)}$ and $x_{Q(\ve v)}$ are relatively prime, we 
obtain that the leading term of   
$ x_{J(\ve v)}\cdot \overline{g}$ divides $y$. This finishes the proof.
\qed


\section{The proof of Theorem \ref{main1}}

Let $0\leq d\leq (q-1)n$. We define the {\em complete $d$-uniform
tuple system}   $\cU(n,d,q)$ as follows:
$$
\cU(n,d,q):=\{\ve v =(v_1,\ldots,v_n)\in (q)^n:~ \sum_{i=1}^n v_i=d\}.
$$

The following result of the authors
from \cite{HR} gives the standard monomials  for the ideal of  
$\cU(n,d,2)$.

\begin{thm} \label{hr} Suppose that $0\leq d\leq n$, and set
$k={\rm min}\{d,n-d\}$. 
Let $\prec$ be an arbitrary term order with $x_n\prec \ldots \prec x_1$.
Then the set of standard monomials of
$\cU(n,d,2)\subset (2)^n $ is
$$
\{x_U:~ U=\{u_1<\cdots <u_\ell \}\mbox{,
  where }\ell \leq k\mbox{ and }u_i\geq 2i\mbox{ for  }1\leq i\leq \ell \}.
$$ \qed
\end{thm}

The sets $U$ appearing in the theorem are essentially the {\em ballot 
sequences} (see \cite{M} or \cite{R}): the characteristic vector of $U$, when viewed as a sequence, has 
at least as many zeros as ones in any initial segment.   

We shall use the approach of \cite{HR} to obtain an upper bound for the 
low degree standard monomials of   $I(\cU(n,d,q))$. 
First we set 
$$
\cB=\cB(n,q)= \{\monom xv:~  \ve v\in(q)^n, ~~ |\{i\leq 2t-1:v_i=q-1\}|\leq
t-1 \mbox{ for all  } t \}.
$$

Next we recall the definition of $\cH(t)$ from \cite{HR}, where it was used
in the description of the leading monomials for $\cU(n,d,2)$. 
Let $t$ be a integer, $0<t\leq n/2$. We define
$\cH(t)$ as the set of those subsets $H=\{s_1<s_2<\cdots <s_t\}$
of $[n]$ for which $t$ is the smallest index $j$ with  $s_j<2j$.
Thus, the elements of $\cH(t)$ are $t$-subsets of $[n]$. We have 
$H\in\cH(t)$ iff $s_1\geq 2,\ldots, s_{t-1}\geq 2t-2$ and $s_t<2t$. 
For the first few values of $t$ it is easy to give $\cH(t)$ explicitly:
we have $\cH (1)=\{\{1\}\}$,
$\cH (2)=\{\{2,3\}\}$, and $\cH(3)=\{\{2,4,5\},\{3,4,5\}\}$.

Now let $0<t\leq n/2$, $0\leq d\leq n$ and $H\in \cH(t)$. Put
$$
H'=H\cup \{2t,2t+1,\ldots ,n\}\subseteq [n].
$$

Let $\cB^c$ stand for the set of monomials in $R$ which are not in 
$\cB$.

\begin{prop}  \label{leadingmonomial}
We have $\cB ^c \subseteq {\rm LM}(I(\cU(n,d,q)))$.
\end{prop}
\proof 
Let $\monom xv \in \cB^c$, with $\ve v=(v_1, \ldots, v_n)$. If there is an $i$ such that 
$v_i\geq q$, then the statement is obvious, $\monom xv$ is a leading monomial 
even for
$I((q)^n)$.  We can therefore assume that $\ve v\in (q)^n $. 
We define now the following tuple $\ve w=(w_1,\ldots, w_n)\in
(2)^n :$ 
\[
w_i:=\left\{ \begin{array}{ll}
1 & \textrm{if $v_i=q-1$} \\
0 & \textrm{if $v_i< q-1$.} 
\end{array} \right.
\]

Let $F=F_{\ve v}$ be the unique subset of $[n]$ such that $\ve w=\ve v_F$, 
where $\ve v_{F}$ stands for the characteristic vector of the set $F$. 
By our assumption on $\monom xv$,  there exists a positive integer $t$
and a $H\in
\cH(t)$ such that $H\subseteq F_{\ve v}$. 
Then, writing 
 $H=\{h_1<\cdots <h_t\}$,
 from the definition of $F_{\ve v}$ we see that
$x_{h_1}^{q-1}\cdots x_ {h_t}^{q-1}$ divides $\monom xv$.  Thus, it suffices 
to prove
that  $x_{h_1}^{q-1}\cdots x_{h_{t-1}}^{q-1}x_{h_t}\in
\mbox{LM}(I(\cU(n,d,q)))$, because then $\monom xv \in
\mbox{LM}(I(\cU(n,d,q)))$ holds as well.

Consider the following polynomial:
$$
f(x_1,\ldots, x_n):=\prod_{i=0}^{(q-1)(t-1)} \left( 
\sum_{h\in H'} x_h-(d-i)\right) \in \Q[x_1,\ldots ,x_n].
$$  
We claim that
$$
f\in I(\cU(n,d,q)).
$$
Indeed, let $\ve u=(u_1,\ldots ,u_n)\in \cU(n,d,q)$ be an arbitrary tuple. 
Then
\begin{equation} \label{sum}
\sum_{h\in H'} u_h=\sum_{i=1}^n u_i -\sum_{j\in [n]\setminus H'} u_j=
d-\sum_{j\in [n]\setminus H'} u_j.
\end{equation}

But clearly $|[n]\setminus H'|=t-1$, therefore 
\begin{equation} \label{sumy}
0\leq \sum_{j\in [n]\setminus H'} u_j\leq (t-1)(q-1).
\end{equation}
Equations (\ref{sum}) and (\ref{sumy}) imply that 
$$
d-(t-1)(q-1)\leq \sum_{h\in H'} u_h\leq d.
$$
This means that there exists an $i$ with $0\leq i\leq (t-1)(q-1)$ such 
that
$d-i=\sum_{h\in H'} u_h$, giving that $f(\ve u)=0$. 

Let us consider the polynomials $f_i(x_i)$ from (\ref{q-poli})
for $1\leq i\leq n$. Clearly $\cG=\{f_1(x_1),\ldots ,f_n(x_n)\}\subseteq
I(\cU(n,d,q))$. We shall examine the normal form $NF(f,\cG )$ of $f$ 
with respect to $\cG$. We have $NF(f,\cG ) \in \Q[x_1,\ldots ,x_n]$. 

The multinomial theorem gives that 
$\frac{((t-1)(q-1)+1)!}{((q-1)!)^{t-1}}\not=0$ is the coefficient of the 
monomial
$y=x_{h_1}^{q-1}\cdots x_{h_{t-1}}^{q-1}x_{h_t}$ in $f$. This monomial is 
not affected by the reduction process with respect to $\cG$, 
as it is not divisible by 
$x_j^q$ for any $j$, it is of top degree $(q-1)(t-1)+1$ in $f$, and 
because reduction with respect to $\cG$ strictly decreases the degree
(the leading monomial of $f_i$ is the only top degree monomial of $f_i$). 
These imply that $y$ occurs among the monomials of $NF(f,\cG )$ as well. 

In fact, any monomial $y'$ in $NF(f, \cG )$ has total degree at most $(q-1)(t-1)+1$, 
it has degree at most $q-1$ in any of the variables $x_j$. Moreover, it is 
composed of the variables $x_h$, for $h\in H'$. Among these monomials 
$y$ is obviously the largest one with respect to $\prec$.   
This implies that $y$ is the leading monomial of $NF(f, \cG )$:  
\begin{equation} \label{egyen} 
\mbox{lm}(NF(f,\cG ))=x_{h_1}^{q-1}\cdots x_{h_{t-1}}^{q-1}x_{h_t}.
\end{equation}

Moreover,  
$f\in I(\cU(n,d,q))$ and $\cG \subseteq I(\cU(n,d,q))$ imply
that  $NF(f, \cG )\in  I(\cU(n,d,q))$.  This fact and (\ref{egyen}) show that 
$$
x_{h_1}^{q-1}\cdots 
x_{h_{t-1}}^{q-1}x_{h_t}=
\mbox{lm}(NF(f, \cG ))
\in \mbox{LM}(I(\cU(n,d,q))).
$$
\qed

\begin{cor} \label{main}
We have 
${\rm SM}(I(\cU(n,d,q)))\subseteq \cB$. \qed
\end{cor}

For an integer $0\leq i\leq n$ we set 
$$ X_i=X_i(n,q):=\{\monom xu:~ \ve u=(u_1\ldots,u_n)\in(q)^n \mbox{ and }|\{j:~
u_j=q-1\}|=i\},$$
and similarly
$$ X_{\leq i}=X_{\leq i}(n,q):=
\{\monom xu:~ \ve u=(u_1\ldots,u_n)\in(q)^n \mbox{ and }|\{j:~
u_j=q-1\}|\leq i\}.$$

\begin{prop} \label{becsles}
Let  
$\cV\subseteq \cU(n,d,q)$ be a $d$-uniform family such that 
 $\mbox{\rm Sh}(\cV)\subseteq {[n] \choose {\leq s}}$. Then
$$
|\cV|\leq |{\rm SM}(I(\cU(n,d,q)))\cap X_{\leq s}|.
$$
\end{prop}

\proof Since $\cV\subseteq \cU(n,d,q)$, we 
have $\mbox{SM}(I(\cV))\subseteq \mbox{SM}(I(\cU(n,d,q)))$. 
On the other hand, $\mbox{Sh}(\cV)\subseteq {[n] \choose {\leq s}}$ implies 
by Proposition \ref{alap} that 
$\mbox{SM}(I(\cV))\subseteq X_{\leq s}$.

We obtain that $\mbox{SM}(I(\cV))\subseteq 
\mbox{SM}(I(\cU(n,d,q)))\cap X_{\leq s}$. The statement now follows, 
since by (\ref{egyenlo}) we have
$|\cV|=|\mbox{SM}(I(\cV))|$.
\qed

\begin{lemma} \label{szam} 
Suppose that  $0\leq i\leq n/2$. Then
$$
|\cB\cap\ X_i|=(q-1)^{n-i}\left( {n \choose i}-{n\choose  i-1}\right) .
$$
\end{lemma}
\proof We set
$$
\cW(q,i):=\{\ve w\in(q)^n:~ \ve x^{\ve w} \in \cB\cap X_i \}.
$$
Obviously we have $|\cB\cap X_i|=|\cW(q,i)|$. The elements of $\cW(q,i)$
are $q$-ary analogs of ballot sequences: in each initial segment they have
at least as many components with value less than $q-1$ as components with
value exactly $q-1$; moreover, the total number of components of the latter
type is $i$.

Consider now the following map $F$ from $(q)^n$ to $(2)^n$:
\[
F(\ve v)_i:=\left\{ \begin{array}{ll}
1 & \textrm{if $v_i=q-1$} \\
0 & \textrm{if $v_i< q-1$.} 
\end{array} \right.
\]
We observe that $G:=F|_{\cW(q,i)}: \cW (q,i)\to \cW(2,i)$ is onto, and 
that $|G^{-1}(\ve u)|=(q-1)^{n-i}$ for each $\ve u\in \cW(2,i)$.
 
The determination of $|\cW(2,i)|$ is the classical problem of counting  
ballot sequences. It is well--known 
(see Theorem 1.1 in \cite{M} or \cite{R}) that
$$
|\cW(2,i)|={n \choose i}-{n \choose {i-1}},
$$
hence $|\cB\cap X_i|=|\cW(q,i)|=$
$$
\sum_{\ve u\in \cW(2,i)} |G^{-1}(\ve u)|=
(q-1)^{n-i}\cdot |\cW(2,i)|
=(q-1)^{n-i}\left( {n \choose i}-{n \choose i-1}\right) .
$$
\qed

To conclude the proof of Theorem \ref{main1}, it suffices to verify 
that if $s$ is an integer, $0\leq s\leq n/2$,  then  
$$
|{\rm SM}(I(\cU(n,d,q)))
\cap X_{\leq s} |\leq
\sum_{i=0}^s (q-1)^{n-i}\left ({n \choose i}-{n\choose  i-1}\right) .
$$

Indeed, we have $\mbox{SM}(I(\cU(n,d,q)))\subseteq \cB$ by Corollary 
\ref{main}, hence 
$$
\mbox{SM}(I(\cU(n,d,q)))\cap X_{\leq s}
\subseteq \cB\cap X_{\leq s}. $$
Therefore it is enough to see that
$$
|\cB\cap X_{\leq s}|\leq \sum_{i=0}^s (q-1)^{n-i}\left( 
{n \choose i}-{n\choose  i-1}\right) .
$$
But this follows at once from Lemma \ref{szam} and the 
disjoint union decomposition below 
$$
\cB\cap X_{\leq s}
=\cup_{i=0}^s( \cB\cap X_i).
$$
This concludes the proof of Theorem \ref{main1}.
\qed

\section{Standard monomials for Hamming spheres}

Our main objective here is to prove Theorem \ref{main2}. To this end it will 
be useful to consider the $q$-ary Hamming spheres: let $0\leq d\leq n$, and 
$$
\cV(n,d,q):=\{\ve v =(v_1,\ldots,v_n)\in (q)^n:~ |\{i\in[n]:~ v_i\neq
0\}|=d\}.
$$
We shall first describe the standard monomials for $I(\cV(n,d,q))$. 
This will extend the corresponding result of  \cite{ARS} to a multivalued
setting.

From $\cU(n,d,2)=\cV(n,d,2)$ and Theorem \ref{hr}
the next statement is immediate.  

\begin{cor} \label{kicsi}
If $0\leq s\leq {\rm min}\{d,n-d\}$,
then the standard monomials of $\cV(n,d,2)$ of degree 
at most $s$ are exactly the standard monomials of $\cV(n,s,2)$. \qed
\end{cor}

By exploiting the relation $\cV(n,d,q)=\cV(n,d,2)^q$
we can now explicitly describe the normal set of $I(\cV(n,d,q))$. 

\begin{cor} \label{unif}
Let $\ve u=(u_1, \ldots, u_n)\in (q)^n$, and set $c:=|J(\ve u)|$.
We have  $\monom xu \in
{\rm SM}(I(\cV(n,d,q)))$ iff the following two conditions are satisfied:
\\
a)\  $c\leq d$ and $|Q(\ve u)|\leq {\rm min}(d-c,n-d)$. \\
b)\ If we write $Q(\ve u)\cup Z(\ve u)$ in the form 
$\{j_1<\ldots <j_{n-c}\}$,  
and if $Q(\ve u)=\{j_{m_1}<\ldots <j_{m_\ell}\}$, then 
$m_{i}\geq 2i$ holds for 
every $1\leq i \leq \ell$. 
\end{cor}

\proof We have  $\cF^q=\cV(n,d,q)$, 
where $\cF:={[n] \choose d}$. For $J\subseteq [n]$ we have
$$
\cF_J=\{F\subseteq [n]:~ |F|=d, \mbox{ and } F\supseteq J\},
$$
hence $\cF_J\neq \emptyset$ iff $|J|\leq d$. From Theorem \ref{monom} we 
obtain that
$$
\mbox{SM}(I(\cV(n,d,q)))=\{\monom xu :~~\ve u\in (q)^n,~|J(\ve u)|\leq d  
\mbox{ and }
x_{Q(\ve u)}\in \mbox{SM}(I(\cF_{J(\ve u)}))\}.
$$
The standard monomials of $\cF_{J(\ve u)}$ are the same as 
the standard monomials 
of the family of all $d-c$-subsets of the set $Q(\ve u)\cup Z(\ve u)$.  
Theorem  \ref{hr} gives now the statement.
\qed

The following upper bound is a consequence of the description of the normal
set $\mbox{SM}(I(\cV(n,d,q)))$ given in Corollary \ref{unif}.

\begin{lemma} \label{becsles3}
Let $0\leq s,d\leq n$, $n\geq 3$, $q\ge 3$ be integers. Suppose that 
$0\leq s+d\leq n$. Then
$$
|{\rm SM}(I(\cV(n,d,q)))\cap X_{\leq s}| 
\leq  {n \choose s} \sum_{i=0}^d (q-2)^i {n-s \choose i}.
$$
\end{lemma}

\proof For $0\leq i\leq d$ we set
$$
\cM_i:=\{\monom xu \in\mbox{SM}(I(\cV(n,d,q))):~ |J(\ve u)|=i\}.
$$
From Corollary \ref{unif} it is easy to verify that 
\begin{equation} \label{Mi}
|\cM_i|={n\choose i}(q-2)^i{n-i \choose {d-i}}={n\choose d}{d\choose i}(q-2)^i.
\end{equation}
From Corollary \ref{kicsi} we know that if $0\leq \ell \leq m$ and
$0\leq s\leq \min(\ell,m-\ell)$, then 
\begin{equation} \label{smaller}
|\{ y\in \mbox{SM}(I(\cV(m,\ell ,2))), \mbox{~deg} \, y\leq s\} |={m\choose s}.
\end{equation}
For $0\leq i\leq d$ we now set
$$
\cN_i:=\cM_i\cap  X_{\leq s}.
$$
Using Corollary \ref{unif}, formulae (\ref{smaller}) and (\ref{Mi}) we obtain 
that   
\begin{equation} \label{metszet}
|\cN_i| =\left\{ \begin{array}{ll}
{n\choose i}{n-i\choose s}(q-2)^i={n\choose s}{n-s \choose i}(q-2)^i &
\textrm{if $s\leq \min(d-i,n-d)$} \\
{n\choose i}{n-i\choose d-i}(q-2)^i={n\choose d}{d \choose i}(q-2)^i & \textrm{otherwise.} 
\end{array} \right.
\end{equation}
We have  
$$
\mbox{Sm}(I(\cV(n,d,q)))\cap  X_{\leq s} =\cup_{i=0}^d \cN_i,
$$
hence it suffices to give an upper bound 
for $\sum_{i=0}^d |\cN_i|$.

\medskip

\noindent
{\bf Claim.} For $0\leq i\leq d$ we have 
$$
|\cN_i|\leq {n\choose s}{n-s \choose i}(q-2)^i.
$$
\proof
First suppose that $d-i\leq \frac{n-i}{2}$. Then $\min(d-i,n-i-(d-i))=d-i$. 
If $s\leq d-i$, then (\ref{metszet}) gives that 
$|\cN_i|= {n\choose s}{n-s \choose
  i}(q-2)^i$. But if $s>d-i$, then using that $s\leq n-d=n-i-(d-i)$, 
we get ${n-i
  \choose d-i}\leq {n-i \choose s}$, implying that 
$$
|\cN_i|={n-i \choose d-i}{n\choose i}(q-2)^i\leq {n-i\choose s}{n\choose i}(q-2)^i={n\choose s}{n-s \choose i}(q-2)^i.
$$

Suppose now that  $d-i> \frac{n-i}{2}$. Then
$\min(d-i,n-i-(d-i))=n-d$. Since $s\leq n-d$, equation (\ref{metszet}) implies that 
$$
|\cN_i|={n\choose s}{n-s \choose i}(q-2)^i,
$$
and this gives the claim.
\qed

We conclude that
$$
\sum_{i=0}^d |\cN_i|\leq \sum_{i=0}^d 
{n\choose s}{n-s \choose i}(q-2)^i={n\choose s} 
\sum_{i=0}^d {n-s \choose i}(q-2)^i. 
$$
This finishes the proof of the Lemma. \qed

We are prepared  now to prove Theorem \ref{main2}. \\

{\bf Proof of  Theorem \ref{main2}:} As the result is known to hold for 
$q=2$, we can assume that $q>2$. Since 
$\cV\subseteq \cV(n,d,q)$, we have also 
$$\mbox{SM}(I(\cV))\subseteq \mbox{SM}(I(\cV(n,d,q))).$$ 
On the other hand, $\cV$ does not shatter sets of size $s+1$, hence by 
Proposition \ref{alap} we obtain that 
$$
\mbox{SM}(I(\cV))\subseteq  X_{\leq s}.
$$
Using  Lemma \ref{becsles3} we obtain
$$
|\cV|=|\mbox{SM}(I( \cV))|\leq |\mbox{SM}(I(\cV(n,d,q)))\cap
 X_{\leq s} |
\leq
$$
$$
\leq {n\choose s} \sum_{i=0}^d {n-s \choose i}(q-2)^i.
$$
This finishes the proof of the theorem.  \qed

\medskip

\section{Concluding remarks}

\noindent
1. Most of our results are also valid over over fields other than $\Q$.
We call a field $\F$ {\em large}, if
the characteristic of $\F$ is 0 or at least $q$. If $\F$ is a
large field, then we can consider $(q)^n$ as a subset of $\F^n$ in a natural
way. The statements in Sections 2 and 4 and the proofs we have given there 
are all valid over arbitrary large fields.

\medskip

\noindent
2. We developed a Gr\"obner basis approach to study shattering in a
multivalued setting. We remark here that the main result of Alon \cite{A}
also has a quite natural and simple proof in the framework of standard 
monomials. 

Alon's Theorem states, that for every tuple system $\cV\subseteq (q)^n$ 
there exists a downward closed tuple system $\cW\subseteq (q)^n$ such that 
$|\cV |=|\cW | $ and for every $S\subseteq [n]$ we have

$$ |\{\ve v|_S:~\ve v\in \cW \}|\leq |\{\ve v|_S:~\ve v\in \cV \}|. $$

In fact, let $\F$ be a large field, and $\prec$ an arbitrary term order 
on the polynomial ring $\F[x_1,\ldots ,x_n]$. We have then 
$\cV\subseteq (q)^n\subseteq \F^n $, and we can consider the set of 
standard monomials   ${\rm SM}(I(\cV))$. One can verify that the set of 
exponent vectors 
$$ \cW=\{\ve u\in (q)^n:~\monom xu\in  {\rm SM}(I(\cV))  \} $$
will meet the requirements of Alon's Theorem\footnote{Alon's
Theorem is formulated in a slightly more general setting, where $\cV$ and 
$\cW$ are subsets of $(q_1)\times (q_2)\times \cdots \times (q_n)$, where the
$q_i$ are positive integers. The 
proof outlined here can be extended without much 
difficulty to the more 
general case.}.  Indeed, it is obvious that $\cW$ is downward closed and 
$|\cV |=|\cW | $. Also, suppose that $S\subseteq [n]$, and  
let  
$$ \cU =\{\ve u\in \cW:~{\rm supp}(\ve u)\subseteq S \}.$$
Using that $\cW$ is downward closed, we see that 
$|\{\ve v|S:~\ve v \in \cW \}|=|\cU|.$
Finally, the set of monomials $\{\monom  xv :~\ve v \in \cU \}$ is linearly
independent on $\cV$, and therefore on 
$\{\ve v|_S:~\ve v\in \cV \}$ as well.

\medskip

\noindent
3. To complement Theorem \ref{main1}, we give here a simple lower bound for 
the size of a $d$-uniform
tuple system $\cV$, which does not shatter an $(s+1)$-element set. 
We start with the following set of tuples (which shows that the
Karpovsky--Milman Theorem is sharp):
$$
\cW(n,s,q):=\{\ve u=(u_1,\ldots,u_n)\in (q)^n:~ |\{i:~ u_i=q-1\}|\leq s\}.
$$
It is immediate that
$$
|\cW(n,s,q)|=\sum_{i=0}^s (q-1)^{n-i}{n \choose i}.
$$
The union
$$
\cW(n,s,q)=\cup_{j=0}^{(q-1)n} (\cW(n,s,q)\cap \cU(n,j,q))
$$
is disjoint. This implies the existence of a $d$ such that
$0\leq d\leq (q-1)n$  and
$$
|\cW(n,s,q)\cap \cU(n,d,q)|\geq \frac{\sum_{i=0}^s (q-1)^{n-i}{n
\choose i}}{(q-1)n+1}.
$$
Clearly $\cX:=\cW(n,s,q)\cap \cU(n,d,q)$ is $d$-uniform and
$\mbox{sh}(\cX)\subseteq {[n] \choose {\leq s}}$.

\medskip

\noindent
4. 
We can easily see that if $s>\lceil \frac{d}{q-1}\rceil$ and
$\cV\subseteq \cU(n,d,q)$ is an arbitrary $d$-uniform tuple system,
then $S\notin \mbox{Sh}(\cV)$, whenever  $S\subseteq [n]$, $|S|=s$. 
For contradiction, suppose that there
exists an $S\in {[n] \choose s}$ such that $S\in \mbox{Sh}(\cV)$.
Define $\ve v =(v_1, \ldots ,v_n)$ as
$$
v_j:=\left\{ \begin{array}{ll}
q-1, & \textrm{if $j\in S$} \\
0, & \textrm{if $j\in [n]\setminus S$.}
\end{array} \right.
$$
Then
$$  
\sum_{i\in S} v_i=s(q-1)>\lceil \frac{d}{q-1}\rceil (q-1)\geq d .
$$
  
From $S\in \mbox{Sh}(\cV)$ we have that
there exists  a $\ve u\in \cV$ such that $\ve u\mid_S=\ve v\mid_S$. Then we
have
$$  
\sum_{i\in S} u_i=\sum_{i\in S} v_i > d=\sum_{i=1}^n u_i,
$$
a contradiction.

\medskip

\noindent
5. The bound of Theorem \ref{main2} is sharp in the case
$n=s+d$, $d\leq n/2$, as witnessed by $\cV:=\cV(n,d,q)$. 
The result is not sharp
for $q>2$ and $s+d<n$, as in this case the last inequality
in the proof is strict.

We remark, that by $d\leq n-s$ we also have the simpler inequality
$$|\cV|
\leq {n\choose s} \sum_{i=0}^d {n-s \choose i}(q-2)^i\leq {n\choose s}
\sum_{i=0}^{n-s} {n-s \choose i}(q-2)^i= {n\choose s}(q-1)^{n-s} .
$$
For $q=2$ this simpler inequality
gives back essentially  the Frankl--Pach bound.

\medskip

\noindent
6. 
For a recent partial improvement of the Frankl--Pach bound we refer to 
Mubayi and Zhao \cite{MZ}.
Shattering and related notions have many important applications in
mathematics and computer science. The interested reader is referred to
Babai and Frankl \cite{BF}, F\"uredi and Pach \cite{FPa},
and Vapnik \cite{V} for more details.

\end{document}